\documentclass[12pt]{article}
\usepackage{amsmath}
\usepackage{amsthm}
\usepackage{amsfonts}
\usepackage{amssymb}

\input epsf

\newtheorem{thm}{Theorem}[section]
\newtheorem{lemma}[thm]{Lemma}

\newtheorem{corollary}[thm]{Corollary}
\newtheorem{proposition}[thm]{Proposition}

\newtheorem{definition}{Definition}

\newcommand{\C}{{\mathbb C}}

\newcommand{\tr}{\text{tr}\;}

\def\hs{\hfill \square}

\begin{document}

\title{The two-eigenvalue problem and density of Jones representation of
braid groups}
\author{Michael H. Freedman$^{\dag}$, Michael J. Larsen$^{\ddag}$,
and Zhenghan Wang$^{\ddag}$}

\maketitle $\dag$ {\it Microsoft Research, One Microsoft Way,
michaelf@microsoft.com}

$\ddag$ {\it Indiana Univ., larsen@math.indiana.edu and
zhewang@indiana.edu}

\vskip.8in

\centerline{\bf Contents} \vskip.3in

Introduction

1.  The two-eigenvalue problem

2.  Hecke algebra representations of braid groups

3.  Duality of Jones-Wenzl representations

4.  Closed images of Jones-Wenzl sectors

5.  Distribution of evaluations of Jones polynomials

6.  Fibonacci representations

\vskip.2in

\centerline{\bf Introduction}

In 1983 V.~Jones discovered a new family of representations $\rho$
of the braid groups.  They emerged from the study of operator
algebras (type $\Pi_1$ factors) and unlike earlier braid
representations had no naive homological interpretation.  Almost
immediately he found that the trace or $\lq\lq$Markov" property of
$\rho$ allowed new link invariants to be defined and this ushered
in the era of quantum topology.  There has been an explosion of
link and 3-manifold invariants with beautiful inter-relations,
asymptotic formulae, and enchanting connections to mathematical
physics: Chern-Simons theory and 2-dimensional statistical
mechanics.
 While many sought to bend Jones' theory toward classical topological
objectives, we have found
 that the relation between the Jones polynomial
 and physics allows potentially realistic models of quantum computation
to be
 created [FKW][FLW][FKLW][F].  Unitarity, a hidden
locality, and density of the Jones representation are central to
computational applications.  With this application in mind, we
have returned to some of Jones' earliest questions about these
representations and the distributions of his invariants.  A few
concise answers are stated here in the introduction. Question 9 of
Jones in [J2] asked for the closed images of the irreducible
components of his representation.  We answer Jones' question, and
also identified the closed images for the general $SU(N)$ case
completely.

A salient feature of Jones representation
 is the {\it two-eigenvalue
property:} the image of each braid generator has only two distinct
eigenvalues $\{-1,q\}$. This is obvious from the quadratic Hecke
relation $(\sigma_{i}+1)(\sigma_{i}-q)=0$.  This two-eigenvalue
property plays a key role in the following theorem:

\vskip.1in

\noindent {\bf Theorem 0.1.} {\it Fix an integer $r\geq 5, r\neq
6, 10, n\geq 3$ or $r=10, n\geq 5$.  Let
$$\rho_n^{(2,r)}=\oplus_{\lambda\in
\wedge_{n}^{(2,r)}}\rho_{\lambda}^{(2,r)}:B_n\rightarrow
\prod_{\lambda\in \wedge_{n}^{(2,r)}}U(\lambda)$$ be the unitary
Jones representation of the $n$-strand braid group $B_n$. Then the
closed image $\overline{\rho_n^{(2,r)}(B_n)}$ contains
$\prod_{\lambda\in \wedge_{n}^{(2,r)}}SU(\lambda)$. }

\vskip.1in

Our original motivation for studying Jones representation is for
quantum computation.  The special case $r=5$ has already been used
to show that the $SU(2)$ Witten-Chern-Simons modular functor at
the fifth root of unity is universal for quantum computation
[FLW]. Combining that paper with the above result, we conclude
that the $SU(2)$ Witten-Chern-Simons modular functor at an $r$-th
root of unity is universal for quantum computation if $r\neq
3,4,6$.

Jones was also concerned with the range of values his invariants
assumed and their statistical properties.  For this we must understand
the topology and measure theory of the image $\Gamma$ of $\rho$, since
the Jones polynomial is obtained by tracing them.

There are three levels of detail in the discussion of a finitely
generated group (or semi-group) $\Gamma$ approximating a Lie group
$G$.  First is density and the rate at which density is achieved.
From [Ki][So] [NC], we extract:

\vskip.1in

\noindent{\bf Theorem 5.6.} {\it  Let $X$ be a set closed under
inverse in a compact semisimple Lie group $G$ (with Killing
metrics) such that the group closure $\langle X \rangle$ is dense
in $G$. Let $X_l$ be the words of length $\leq l$ in $X$, then
$X_l$ is an $\epsilon$-net in $G$ for $l={\mathcal
O}(log^2(\frac{1}{\epsilon})),$ i.e.,  for all $g \in G$,
$\text{dist}(g,X_l)<\epsilon$. }

\vskip.1in Conjecturally the theorem should still hold for
$l={\mathcal O}(log{(\frac{1}{\epsilon})})$ and there are some
number theoretically special generating sets of $SU(2)$ [GJS] for
which such an estimate for $l$ can in fact be obtained.  Such
results now translate into topological statements:

\vskip.1in

\noindent{\bf Corollary 5.7.} {\it Given a $\lq\lq$conceivable"
value $v$ for the evaluation of Jones polynomial of $\hat{b}$ at a
root of unity, i.e., one that lies in the computed support of the
limiting distribution for $b\in B_n$, the $n$-string braids, to
approximate $v$ by $v'$, $||v-v'||<\epsilon$, it is sufficient to
consider braids $b_l'\in B_n$ of length $l={\mathcal
O}(log^2(\frac{1}{\epsilon}))$ with Jones evaluations $b_l'=v'$,
$||v-v'||< \epsilon$. }

 The second level is uniformity in measure: if
$\Gamma=<\gamma_1, \cdots, \gamma_m>$, i.e., $\Gamma$ is generated
as a semi-group by $\gamma_1, \cdots, \gamma_m$, $\;$ let $W_l$ be
the set of unreduced words of length=$l$ and $\mu_l$ be the
equally weighted atomic measure on $W_l$ (mass $m^{-l}$ on each
word in $W_l$),  it is known that density implies uniformity in
measure [Bh],  $\mu_l\rightarrow \text{Haar}(G)$ in the weak-*
topology (i.e., when integrated against continuous functions.)
Third is the rate of convergence of measures, which is also
addressed in [Bh].

Returning to the Jones polynomial evaluations which are weighted
traces of dense representations, we can determine the statistics.
Recall $n$ is the number of strands, and $l$ is the length of a
braid.   One may consider the double limit when $l$, and later $n$
are taken to infinity. In this case,
 if $r$ is a fixed integer $r\ge 5$, $r\neq 6$,
 the distribution of evaluations at $e^{\frac{\pm 2\pi i}{r}}$
of the Jones polynomial
of a ``random'' link with $n$ strands tends to a fixed Gaussian.  The
variance
of this Gaussian depends on $r$ and grows like $r^3$ as $r\to\infty$.

\vskip.1in Our density result follows from the solution of a
general two-eigenvalue problem: Let $G$ be a compact Lie group,
and $V$ a faithful, irreducible, unitary representation of $G$.
The pair $(G,V)$ is said to have the {\it $k$-eigenvalue property}
if there exists a conjugacy class $[g]$ of $G$ such that

(1) the class $[g]$ generates $G$ topologically;

(2) any element $g \in [g]$ acts on $V$ with exactly $k$ different
eigenvalues such that for each $2\leq r\leq k$, no set of $r$
eigenvalues forms a coset of the multiplicative group $\{1,
\omega, \omega^2, \cdots , \omega^{r-1} \}$, where $\omega$ is a
primitive $r$-th root of unity.

The $k$-eigenvalue problem is to classify all such pairs $(G,V)$.
Note that $G$ is not assumed to be connected.  The problem
naturally divides into two cases according to whether $G$ is or is
not finite modulo its center.  The solution to the first case is
essentially known to the experts and we content ourselves with a
statement at the end of section 1.  The solution to the case that
$G/Z(G)$ has positive dimension is:

\vskip.1in

\noindent{\bf Theorem 1.1} {\it Suppose $(G,V)$ is a pair with the
two-eigenvalue property. Let $G_1$ be the universal covering of
the derived group $[G_0,G_0]$ of the identity component $G_0$ of
$G$.  If $G$ is of positive dimension modulo its center,
 then $V$ is an irreducible $G_1$-module, with highest weight $\varpi$,
and $(G_{1},\varpi)$ is one of the following:

(1) $(SU(l+1), \varpi_{i})$ for some $l\geq 1$,  and $1\leq i\leq
l.$

(2) $(Spin(2l+1), \varpi_{l})$ for some $l\geq 2$.

(3) $(Sp(2l), \varpi_{1})$ for some $l\geq 3$.

(4) $(Spin(2l), \varpi_{i})$ for some $l\geq 4$ and $ i=1,l-1,l,$

\noindent where $\varpi_i$ denotes the $i$-th fundamental
representation.

} \vskip.1in

There is a fairly close analogy between this theorem and
J.~Serre's classification [Se] of inertial monodromy types for
Hodge-Tate modules with only two different weights.  Not only are
the problems formally similar, the solution is identical. However,
it does not seem that either result implies the other.  In the
Hodge-Tate case, one looks for a cocharacter taking two distinct
values on the set of weights of an irreducible representation of a
semisimple group; in our case, one looks for a rational
cocharacter taking two different values (mod $\mathbb Z$) which
are not congruent (mod $\frac{1}{2}{\mathbb Z}$). Our technique
here works  for the 3-eigenvalue problem.

\section{The two-eigenvalue problem}

Let $G$ be a compact Lie group, and $V$ a faithful, irreducible,
unitary representation of $G$. The pair $(G,V)$ is said to have
the {\it two-eigenvalue property} if there exists a conjugacy class
$[g]$
of $G$ such that

(1) the class $[g]$ generates $G$ topologically;

(2) any element $g \in [g]$ acts on $V$ with exactly two
different eigenvalues whose ratio is not $\pm 1$.

\smallskip
Note that $G$ is not assumed to be connected.  The problem
naturally divides into two cases according to whether $G$ is or is not
finite modulo its center.  The solution to the first case is essentially
known to the experts and we content ourselves with a statement at the
end
of this section.  The rest of the section is devoted to the case that
$G/Z(G)$
has positive dimension.

\begin{thm}\label{twoeigenvalue}
Suppose $(G,V)$ is a pair with the two-eigenvalue property.  Let
$G_1$ be the universal covering of the derived group $[G_0,G_0]$
of the identity component $G_0$ of $G$.  If $G$ is
of positive dimension modulo its center,
 then $V$ is an irreducible $G_1$-module, with highest weight $\varpi$,
and $(G_{1},\varpi)$ is one of the following:

(1) $(SU(l+1), \varpi_{i})$ for some $l\geq 1$,  and $1\leq i\leq
l.$

(2) $(Spin(2l+1), \varpi_{l})$ for some $l\geq 2$.

(3) $(Sp(2l), \varpi_{1})$ for some $l\geq 3$.

(4) $(Spin(2l), \varpi_{i})$ for some $l\geq 4$ and $ i=1,l-1,l,$

\noindent where $\varpi_i$ denotes the $i$-th fundamental
representation.

\end{thm}

There is a fairly close analogy between this theorem and
J.~Serre's classification [Se] of inertial monodromy types for
Hodge-Tate modules with only two different weights.  Not only are
the problems formally similar, the solution is identical. However,
it does not seem that either result implies the other.  In the
Hodge-Tate case, one looks for a cocharacter taking two distinct
values on the set of weights of an irreducible representation of a
semisimple group; in our case, one looks for a rational
cocharacter taking two different values (mod $\mathbb Z$) which
are not congruent (mod $\frac{1}{2}{\mathbb Z}$).



We begin with a lemma from linear algebra.

\begin{lemma}
 Suppose $W$ is a vector space with a direct sum decomposition
$W=\oplus_{i=1}^{n} W_{i}$, and $U$ is an operator on
$W$ such that $U:W_{i}\rightarrow W_{i+1}$ $(1\leq i\leq n)$
cyclically. Then any eigenvalue of $U$ multiplied by any $n$th
root of unity is again an eigenvalue of $U$.
\end{lemma}

{\bf Proof:} Choose a basis of $W$ consisting of bases of
$W_{i},i=1,2,\cdots, n$. If $k$ is not a multiple of $n$, then
$\tr U^k=0$ because all diagonal entries of $U^k$ are 0 with
respect to the above basis.  Let $\lambda_1,\ldots,\lambda_N$
denote the eigenvalues of $U$ with multiplicity. For each integer
$m>0$, consider $\tr U^m=\sum {\lambda_{i}}^{m}$.  Let $\omega$ be
an $n$th root of unity.  Then $\sum {(\omega
\lambda_{i})}^{m}=\sum\omega^{m}
{\lambda_{i}}^{m}={\omega}^{m}\sum {\lambda_{i}}^{m}$.    We claim
this sum is equal to $\tr U^m=\sum {\lambda_{i}}^{m}$.  Indeed,
when $m$ is not a multiple of $n$, they are both $0$,  when $m$ is
a multiple of $n$, $\omega^{m}=1$. Recall that the symmetric
polynomials  $\{\sum x_i^m\} $ uniquely determine all the
symmetric polynomials of $x_i$. It follows that
$\prod_{i}(\lambda-\omega
\lambda_{i})=\prod_{i}(\lambda-\lambda_{i})$. Therefore, the set
of the eigenvalues of $T$ is invariant under multiplication by any
$n$th root of unity. $\hs$

 In the two-eigenvalue problem, the generating conjugacy class cannot
lie
in the identity component $G_0$ unless $G$ is connected.  However, the
following lemma allows us to reduce to the connected case:

\begin{lemma}
Given a compact Lie group $G$, and  an irreducible representation
of $G$. If an element $g$ has two eigenvalues under $\rho$ whose
ratio is $\neq \pm 1$, then $g$ is a product of an element in $G_0$
with an element in $Z_{G}(G_0)$,
the centralizer of $G_0$ in $G$.
\end{lemma}

{\bf Proof:} The action of $Ad_g$ defines an automorphism of
$G_0$. By [St]~Theorem~7.5, there exists a maximal torus $T$ of
$G_0$ such that $Ad_g$ fixes $T$ as a set.  Recall any
automorphism of $G_0$ fixing $T$ pointwise is an inner
automorphism by an element in $T$.

To show that $Ad_g$ fixes $T$ pointwise, consider all the
characters $\{\chi\}$ of $\rho$, and the weight space
decomposition $V=\oplus_{\chi \in \chi^{*}(T)}V_{\chi}$.  As
$Ad_g$ fixes $T$ as a set, $\rho(g)$ permutes the weight spaces
$V_{\chi}$ according to the permutation of characters by $Ad_g$.
Suppose the longest permutation cycle of weight spaces by $Ad_g$
has length=$l$.  If $l\geq 3$, then by Lemma 1.2, $\rho(g)$ have
at least $l$ distinct eigenvalues, contrary to hypothesis. If
$l=2$, then by Lemma 1.2, the two possible eigenvalues of
$\rho(g)$ has ratio $-1$.   Therefore, $l=1$, i.e., $\rho(g)$
fixes every weight space $V_{\chi}$. It follows that $Ad_g$ fixes
the maximal torus $T$ of $G_0$ pointwise.  The lemma follows.
$\hs$

\begin{thm}\label{simplicity}
Let $(G,V)$ be a pair with the two eigenvalue property.  If $G$ is
of positive dimension modulo its center, then the
derived group $[G_0,G_0]$ of $G_0$ is a
simple Lie group, and $G=G_0 Z(G)$.
\end{thm}
{\bf Proof: } Let $[g]$ satisfy the two-eigenvalue property.  As the
conjugates of $g$ (topologically) generate $G/G_0$, if the restriction
of $V$ to
$G_0$ had more than one isotypic component, $g$ would permute
these components nontrivially, contrary to Lemma 1.2.  Thus, the
restriction of $V$ to $G_0$ is the tensor product of an
irreducible representation $V_0$ and a trivial representation
$V^0$.  By Lemma 1.3, $g=g_0 z$, where $g_0\in G_0$ and $z$
centralizes $G_0$.  By Schur's Lemma, $\rho(z)=1\otimes B$, while
$\rho(g_0)=A\otimes 1$.  The two-eigenvalue property implies that
either $A$ or $B$ is scalar.  Since
$[g]$ generates a dense subgroup of $G$, the same is true of $[g_0]$ and
$G_0$.  As $V$ is a faithful representation, $A$ cannot be scalar,
so $B$ must be.  Thus, $(G_0,V_0)$ satisfies the two-eigenvalue
property with generating class $[g_0]$. Moreover, $V^0$ must be
one-dimensional since otherwise $V$ would be a reducible
representation of $G$.

Let $G_1$ denote the universal cover of $[G_0,G_0]$.  Let $g_1\in
G_1$ denote an element whose image in $[G_0,G_0]$ lies in the
coset $g_0 Z(G_0)$.  The pull-back $V_1$ of $V_0$ to $G_1$ is
again irreducible, and the image of $g_1$ has two eigenvalues with
the same ratio as the original image of $g_0$. Moreover, $[g_1]$
generates
a dense subgroup of $G_1$
since no proper closed subgroup of $G_1$ can generate $G_0$ modulo
$Z(G_0)$.  It follows that $(G_1,V_1)$ satisfies the
two-eigenvalue property.

If $G_1$ were not simple, it would factor as $G_2\times G_3$, and
$V_1$ would factor as an external tensor product of
representations $V_2$ and $V_3$. Writing $\rho(g_1)=A\otimes B$,
we see that $A$ or $B$ must be a scalar. Thus $[g_1]$ cannot
generate a dense subgroup of the product. We
conclude that $G_1$, and therefore $[G_0,G_0]$, must be simple. $\hs$

\begin{thm} Let $G$ be a connected, simply connected compact simple Lie
group and $V$ an irreducible representation of $G$ satisfying the
two-eigenvalue property.  Let $\varpi$ denote the highest weight
of $V$.  Then $(G,\varpi)$ is one of the following:
\smallskip
\item{(1)} $(SU(r+1),\varpi_i)$ for some $r\ge 1$ and $1\le i\le r$.
\item{(2)} $(Spin (2r+1), \varpi_r)$ for some $r\ge 2$.
\item{(3)} $(Sp(2r), \varpi_1)$ for some $r\ge 3$.
\item{(4)} $(Spin (2r), \varpi_i)$ for some $r\ge 4$ and
$i\in\{1,r-1,r\}$.
\end{thm}

 In other words $G$ is classical and $V$ is minuscule.

\vskip.1in

 {\bf Proof: }  Fix a maximal torus $T$ of $G$.  As the
conjugates of $T$ cover $G$, there exists $g\in T$ satisfying the
two-eigenvalue property. There is a natural identification of $T$
with the quotient $W/X_*(T)$, where $W=X_*(T)\otimes{\bf R}$ is
the universal covering space of $T$, and where we identify ${\bf
R}/{\bf Z}$ with the set of complex numbers of norm 1. Let $\tilde g$
denote
an element of $W$ mapping to $g$.  The
two-eigenvalue condition means that the values $\chi(\tilde g)$,
as $\chi$ ranges over the characters of $V$, lie in exactly two
cosets of ${\bf Z}$ which do not differ by a half-integer.

Let $\alpha$ denote the highest short root of $G$ and
$\varpi,\varpi-\alpha, \ldots, \varpi-k\alpha$ a string of weights
of $V$. If $k\ge 2$, then $\alpha(\tilde g)$ must be an integer.
As the set of weights is invariant under the Weyl group, all short
roots of $G$ lie in the Weyl-orbit of $\alpha$, and as the short
roots span the root lattice, this would imply that all
$\chi(\tilde g)$ lie in a single coset, contrary to hypothesis. It
follows that $k=1$, or equivalently,
$$\sum_{i=1}^{r}a_ib_i\cdot
\frac{\alpha_i^2}{\alpha^2}=1,$$
where
$$\varpi=a_1\varpi_1+\cdots+a_r\varpi_r, \
\alpha=b_1\alpha_1+\cdots+b_r\alpha_r.$$
Indeed, in the notation of [Hu],
$$1=\langle\varpi, \alpha\rangle=2\frac{\varpi\cdot \alpha}{\alpha^2}=
2\sum_{i,j}a_ib_j\frac{\varpi_i \cdot
\alpha_j}{\alpha^2}=\sum_{i,j}a_ib_j\langle\varpi_i,
\alpha_j\rangle\frac{\alpha_j^2}{\alpha^2}=
\sum_{i}a_ib_i\frac{\alpha_i^2}{\alpha^2}.$$
 Note that
$\frac{\alpha_i^2}{\alpha^2}\in \{1,2,3\}$. Since all the
coefficients $b_i$ in the representation of the longest short root
as a linear combination of simple roots are $\ge 1$, this implies
that $\varpi$ is a fundamental weight $\varpi_i$ for some $i$ such
that $a_i=b_i=1$, and $\alpha_i$ is a short root. In addition to
the cases listed above, we have the cases $(E_6,\varpi_1)$,
$(E_6,\varpi_6)$,
and $(E_7,\varpi_7)$. We claim that none of these exceptional cases
correspond to actual solutions of the two-eigenvalue problem.

For $E_6$, the two representations in question are dual to one
another, so we consider only the one corresponding to the
highest weight $\varpi_1$. By [MP], the restriction of this
representation to $H=SU(3)\times SU(3)\times SU(3)$ is
$$\sigma\otimes \sigma^*\otimes 1\oplus
1\otimes\sigma\otimes\sigma^* \oplus\sigma^*\otimes
1\otimes\sigma,$$
where $\sigma$ denotes the standard
representation of $SU(3)$.  Since $H$ can be chosen to contain
$T$, we may write $g=(g_1,g_2,g_3)\in H$. The two-eigenvalue
property guarantees that one of the $\sigma(g_i)$ has two
eigenvalues and the other two are scalars.  Without loss of
generality, we assume $\sigma(g_1)$ has eigenvalues $\alpha$ (with
multiplicity 2) and $\alpha^{-2}$, while the scalars for $g_2$ and
$g_3$ are $\beta$ and $\gamma$. The set of eigenvalues is
$$\{\alpha\beta^{-1},\alpha^{-2}\beta^{-1},\beta\gamma^{-1},\gamma\alpha
^{-1},
\gamma\alpha^2\}.$$
Since two pairs of eigenvalues have ratio
$\alpha^3$, either $\alpha\beta^{-1}=\gamma\alpha^2$ or
$\alpha^3=1$.  In the first case, $\alpha\beta\gamma=1$, and since
$\beta^3=\gamma^3=1$, this implies $\alpha^3=1$.  We conclude that
the eigenvalues are $\alpha/\beta$, $\beta/\gamma$, and
$\gamma/\alpha$, all cube roots of unity.  Since they multiply to
1, all are the same or all are different, contrary to hypothesis.

For $E_7$, we restrict to $SU(2)\times SU(4)\times SU(4)$ and
obtain
$$1\otimes\sigma\otimes\sigma\oplus
1\otimes\sigma^*\otimes\sigma^* \oplus\tau\otimes 1\otimes
S^2\sigma\oplus\tau\otimes S^2\sigma\otimes 1,$$
where $\sigma$
and $\tau$ are the standard representations of $SU(4)$ and $SU(2)$
respectively.  Writing $g=(g_1,g_2,g_3)$, we conclude that
$\sigma(g_2)$ and $\sigma(g_3)$ are scalars $\beta$ and $\gamma$,
while $\tau(g_1)$ has eigenvalues $\alpha^{\pm1}$.  Thus, the set
of eigenvalues is
$$\{\beta\gamma,\beta^{-1}\gamma^{-1},\alpha\gamma^2,\alpha^{-1}\gamma^2
,
\alpha\beta^2, \alpha^{-1}\beta^2\}.$$
Note that $\gamma^2=\beta^2=\pm1$ since $\beta$ and $\gamma$
determine unimodular scalar
$4\times 4$ matrices.  If $\alpha^2=1$, then all the
eigenvalues are the same up to sign, contrary to hypothesis.  If
not the squares of eigenvalues are $1$, $\alpha^2$, and
$\alpha^{-2}$, so $\alpha^2=-1$.  But this implies that two
eigenvalues have ratio $-1$, contrary to hypothesis. $\hs$

Now we state the solution to the two-eigenvalue problem for finite
groups.  Our list is based on [Za] and depends on the
classification of finite simple groups.  The cases $m\ge 5$ are
classical [Bl].

\begin{thm}\label{finite}
Suppose $(G,V,[g])$ has the two-eigenvalue property, and $G/Z(G)$
is finite. Then $g^m\in Z(G)$ for some $m\in \{3,4,5\}$, and
$G=H\cdot Z(G)$ for some group $H$ with an element $h\in H$ such
that $h^{-1}g\in Z(G)$. Furthermore, one of the following holds:
\vskip.1in

\noindent (a) $m=5$, $H \cong SL(2,5)$ and
$\dim V=2$;\vskip.1in

\noindent (b) $m=4$, $G$ contains a normal subgroup $E$ such that
$E/Z(E)$ is of exponent 2 and of order $2^{2k}$,
$\dim V=2^k$, $V|_E$ is irreducible and $H/E \in \{Sp(2k,2),
U(k,2), O^-(2k,2) \;\text{with}\; k>2, S_{2k+1}, S_{2k+2}\}$;

\vskip.1in \noindent (c) $m=3$ and one of the following holds:

(1)  $H\cong Sp(2n,3), n>1$ and
$\dim V=\frac{(3^n-(-1)^n)}{2}$;

(2)  $H\cong PSp(2n,3), n>1$ and
$\dim V=\frac{(3^n+(-1)^n)}{2}$;

(3)  $H\cong SU(n,2)$ and $n$ is a multiple of $3$, or $H\cong
U(n,2)$, $V|_H$ is a Weil representation of $H$ and $\dim V=
\frac{(2^n+2(-1)^n)}{3}$ or $\frac{(2^n-(-1)^n)}{3}$;

(4)  $H\cong \tilde{A_n}$, the two-fold central extension of the
alternating group $A_n$, and $\dim V=2^{\frac{n-3}{2}}$ for n
odd, and $\dim V=2^{\frac{n-2}{2}}$ for n even;

(5)  $G$ contains a normal subgroup $E$ such that $E/Z(E)$ is of
exponent 2 and of order $2^{2k}$, $\dim V=2^k$, $V|_E$ is
irreducible and $H/E$   $  \in \{Sp(2k,2), \linebreak U(k,2),
O^+(2k,2), O^-(2k,2) \;\text{with}\; k>2, A_{2k+1}, A_{2k+2}\}$;

(6)  $G$ contains a normal extraspecial subgroup $E$ of order
$3^{2k}$, $\dim V=3^k$, and $V|_E$ is irreducible, and
$H/E\cong Sp(2k,3)$;

(7)  $H\cong PSp(4,3)$, and $\dim V=6$;

(8)  $H/Z(H)\cong PSU(4,3)$, $|Z(G)|=6$, and $\dim V=6$;

(9)  $H/Z(H)\cong J_2, |Z(G)|=2$, and $\dim V=6$;

(10)  $H/Z(H)\cong Sp(6,2)$, $|Z(G)|=2$, and $\dim V=8$;

(11)  $H/Z(H)\cong O^+(8,2), |Z(G)|=2$, and $\dim V=8$;

(12)  $H/Z(H)\cong G_2(4), |Z(G)|=2$, and $\dim V=12$;

(13)  $H/Z(H)\cong Suz, |Z(G)|=6$, and $\dim V=12$.

(14)  $H\cong Co_1$, and $\dim V=24$;

\end{thm}

\section{Hecke algebra representations of braid groups}

The $n$-strand braid group $B_{n}$ has the well-known
presentation:
$$\begin{array}{cc}
B_n=\{ \sigma_1, \cdots, \sigma_{n-1}| &
\sigma_{i}\sigma_{j}=\sigma_{j}\sigma_{i}\;\; \text{if} \;\;
|i-j|>1  \\
 & \sigma_{i}\sigma_{j}\sigma_{i}=\sigma_{j}\sigma_{i}\sigma_{j}
\;\; \text{if}\;\; |i-j|=1\}.
\end{array} $$
Hecke algebra
representations of the braid groups in the root of unity case
are indexed by two parameters:
a compact Lie group and an integer $l\geq 1$, called the {\it
level} of the theory. The cases of Jones and Wenzl representations
correspond to the special unitary groups $SU(k), k\geq 2$. For
each pair of integers $(k,r)$ with $r\geq k+1$, there is a unitary
representation of the braid groups with level $l=r-k$.  Jones
representations correspond to $SU(2)$, and the general $SU(k)$
theory gives rise to the HOMFLY polynomial.

We describe the Jones-Wenzl representation explicitly, following
[We]. Let $q=e^{\pm \frac{2\pi i}{r}}$, and $[m]$ be the quantum
integer
$\frac{q^{\frac{m}{2}}-q^{-\frac{m}{2}}}{q^{\frac{1}{2}}-q^{-\frac{1}{2}
}}$.
The constant
$[2]=q^{\frac{1}{2}}+q^{-\frac{1}{2}}=2cos\frac{\pi}{r}$ is
ubiquitous in quantum topology.  The Hecke algebra $H_{n}(q)$ of type
$A$ is
the (finite dimensional) complex algebra generated by $e_1, \ldots,
e_{n-1}$ such that
\begin{enumerate}
\item $e_i^2=e_i$,

\item $e_i e_{i+ 1} e_i-[2] ^{-2} e_i=e_{i+1} e_{i} e_{i+1}-[2] ^{-2}
e_{i+1}$,

\item $e_i e_j=e_j e_i$ if $|i-j|\geq 2$.

\end{enumerate}

 A representation $\pi$ of
$H_{n}(q)$ on a Hilbert space is called a {\it $\C^{*}$
representation} if each $\pi(e_{i})$ is self-adjoint.

\begin{lemma} Each $\C^{*}$ representation of the Hecke algebra
$H_{n}(q)$
gives rise to a unitary representation of the braid group $B_{n}$
by the formula:
\begin{equation}
 \rho(\sigma_i)=q-(1+q) \pi(e_i). \label{3.2}
\end{equation}

\end{lemma}

{\bf Proof:}  The defining relations 1---3 of $H_{n}(q)$ imply
that the elements $\rho(\sigma_i)$ satisfy the braid relations.
Writing $e_i$ for $\pi(e_i)$, since
$\rho^{*}(\sigma_{i})=\bar{q}-(1+\bar{q})e^{*}_{i}$,
$$\rho(\sigma_{i})\rho^{*}(\sigma_{i})=
q\bar{q}+(1+q)(1+\bar{q})e_{i}e_{i}^{*}-{\bar{q}}(1+q)e_{i}-
q(1+\bar{q})e_{i}^{*}= 1.$$
Cancellation of the last three terms follows from the facts
$e_{i}^{*}=e_{i}$ and $e_{i}^{2}=e_{i}$. $\hs$

Jones-Wenzl $\C^{*}$ representation of $H_{n}(q)$ are reducible;
their irreducible constituents, referred to as {\it sectors}, are
indexed by Young diagrams. A Young diagram with $n$ boxes is the
diagram of a partition of the integer $n$:
$$\lambda=[\lambda_1,
\ldots, \lambda_k], \lambda_{1}\geq \lambda_{2} \geq \cdots \geq
\lambda_{k}\geq 0,\ \sum_{i=1}^{k}\lambda_{i}=n.$$
Note that
$\lambda$ is allowed to have empty rows.   Given a Young diagram
$\lambda$ with $n$ boxes, a standard tableau of shape $\lambda$ is
an assignment of integers $\{1,2,\cdots, n\}$ into the boxes so
that the entries of each row and column are increasing.

\begin{definition}
 Suppose $t$ is a standard tableau with $n$ boxes, and $m_1$ and
$m_2$ are two entries in $t$. Suppose $m_i$ appears
in row $r_i$ and column $c_i$
of $t$.

(1) Set $d_{t,m_1,m_2}=(c_1 -c_2) -(r_1 -r_2)$.

(2) Set $\alpha_{t,i}=\frac{[d_{t,i,i+1}+1]}{[2][d_{t,i,i+1}]}$
if $[d_{t,i,i+1}]\neq 0$, and
$\beta_{t,i}=\sqrt{\alpha_{t,i}(1-\alpha_{t,i})}$.

(3) A Young diagram $\lambda=[\lambda_1, \cdots \lambda_k],
\lambda_{1}\geq \lambda_{2} \geq \cdots \geq \lambda_{k}\geq 0$
 is $(k,r)$-admissible if  $\lambda_{1}-\lambda_{k} \leq r-k$.

(4) Suppose $t$ is a standard tableau of shape $\lambda$ with
$n$ boxes, let $t^{(i)} (1\leq i\leq n)$ be the standard tableaux
obtained from $t$ by deleting boxes with entries $n,n-1,\cdots,
n-i+1$. A standard tableau $t$ is $(k,r)$-admissible if the shape
of each tableau $t^{(i)}$ is a $(k,r)$-admissible Young diagram.
\end{definition}

The irreducible sectors of the Jones-Wenzl representations of the
Hecke algebras $H_n(q)$ (and hence of the braid groups $B_n$) are
indexed
by the the pair $(k,r)$ and a $(k,r)$-admissible
Young diagram $\lambda$ with $n$ boxes. A $\C^{*}$ representation
${\pi}_{\lambda}^{(k,r)}$ of the Hecke algebra $H_{n}(q)$ can be
constructed as follows: let $V_{\lambda}^{(k,r)}$ be the complex
vector space with basis $\{ \vec{v}_{t}\}$, where $t$ ranges over
$(k,r)$-admissible standard tableaux of shape $\lambda$. Let
$s_{i}(t)$ be the tableau obtained from $t$ by interchanging the
entries $i$ and $i+1$.  If $s_i(t)$ is also $(k,r)$-admissible,
then we define
\begin{equation}
 \pi_{\lambda}^{(k,r)}(e_{i}) (\vec{v}_t)=\alpha_{t,i}\vec{v}_t
+\beta_{t,i}\vec{v}_{s_{i}(t)}. \label{3.0}
\end{equation}
If $s_{i}(t)$ is not $(k,r)$-admissible, set $\beta_{t,i}=0$ in
formula (2). In this case, $\alpha_{t,i}$ is either $0$ or $1$. It
follows that $\pi_{\lambda}^{(k,r)}(e_i)$ (with respect to the
basis $\{\vec{v}_{t} \}$) is a matrix consisting of only $2\times
2$ blocks
\begin{equation}
\left(
\begin{array}{cc}
\alpha_{t,i} & \beta_{t,i} \\
 \beta_{t,i}  & 1-\alpha_{t,i}
\end{array} \right) \label{3.1}
\end{equation}
 and $1 \times 1$ blocks $0$ or $1$.  The identity
$\alpha_{t,i}=\alpha^2_{t,i} + \beta^2_{t,i}$ implies that
(\ref{3.1}) is a projector. So all eigenvalues of $e_{i}$ are
either $0$ or $1$.  We write $\rho_\lambda^{(k,r)}$ for the restriction
of $\pi_\lambda^{(k,r)}$ to $B_n$.  When $n$ and $r$ are fixed,
they may be suppressed.

\begin{definition}
Given a pair of integers $(k,r)$ with $r\geq k+1$.  Let
$\Lambda_{n}^{(k,r)}$ be the set of all $(k,r)$-admissible Young
diagrams with $n$ boxes.  \emph{The Jones-Wenzl representation} of the
braid group $B_n$ is:
$$\rho_{n}^{(k,r)}=\oplus_{\lambda \in
\Lambda_{n}^{(k,r)}}\rho_{\lambda}^{(k,r)}: B_n\rightarrow
\prod_{\lambda \in \Lambda_{n}^{(k,r)}}U(\lambda).$$
Here we write $U(\lambda)$ for the unitary group of
the Hilbert space $V_{\lambda}^{(k,r)}$ with the orthonormal basis
$\{\vec{v}_t\}$.
\end{definition}

\begin{definition}
A $(k,r)$-admissible diagram is of \emph{trivial type} if $\lambda$ is a
row or column or if $k=r-1$.  A $(k,r)$-admissible diagram is a
\emph{hook}
if the second row has exactly one box.  A hook with exactly two rows is
a
\emph{Burau hook}, and the corresponding sector is a \emph{Burau
representation}.
\end{definition}

We note that $\rho_\lambda$ is one-dimensional if and only if $\lambda$
is
of trivial type.

\begin{thm}\label{hook}
Let $h$ be a $(k,r)$-admissible hook with $(b+1)$ rows and
$(a+1)$ columns.

(1) If $a+b<r-1$, then $\rho_{h}^{(k,r)}$ is equivalent up to tensoring
by a character to the $b$th exterior
power of the Burau representation associated to the hook with
$(a+b)$ columns.

(2)  If $a+b=r-1$, then $\rho_{h}^{(k,r)}$ is equivalent up to tensoring
by a character to the $(b-1)$th
exterior power of the Burau representation associated to the hook
with $(a+b-1)$ columns.

\end{thm}

{\bf Proof:} For the first part, we explicitly identify a basis of
$V_{h}$ with
that of the $b$-th exterior power of the Burau representation
$\rho_{\beta}$ associated to the hook $\beta$ with $(a+b)$
columns. The basis of $V_{\beta}$ can be indexed conveniently by
the entry $i$ of the box in the second row.  The set
$$\{ v_{i_2}\wedge
v_{i_3}\wedge \cdots \wedge v_{i_{b+1}} \mid 2\leq i_2<
\cdots  < i_{b+1}\leq a+b+1  \}$$
spans $\wedge^b V_\beta$.  We
identify each element of this basis with the basis element of $V_{h}$
given by the standard tableau whose first column entries are
 $1,i_2,\cdots, i_{b+1}$,
 which we denote $v_{1,i_2,\cdots, i_{b+1}}$.
 Now we just compare
 the action of the braid generator $\sigma_k$ on corresponding
 basis elements: $v_{1,i_2,\cdots, i_{b+1}}$ and $v_{i_2}\wedge
v_{i_3}\wedge \cdots \wedge v_{i_{b+1}}$. For the Burau
representation, we have $\rho_{\beta}(\sigma_{k})(v_i)=
\begin{array}{ccc}
qv_i & \text{if} & i\neq k, k+1
\end{array}.
$
We drop $\rho$ from the notation now.   First we compare two
special cases:
$$\sigma_{k}(v_{1,i_2,\cdots, i_{b+1}})
=\begin{cases}
q &\text{if $k$ and $k+1$ do not appear in $i_2,\ldots,
i_{b+1}$} \\
-1 &\text{if $k$ and $k+1$ both appear in  $i_2,\ldots,
i_{b+1}$} \\
\end{cases}
$$
$$\sigma_{k}(v_{i_2}\wedge \cdots \wedge v_{i_{b+1}})=\begin{cases}
q^b &\text{if $k$ and $k+1$ do not appear in $i_2, \ldots,
i_{b+1}$} \\
 -q^{b-1} &\text{if $k$ and $k+1$ both appear in $i_2, \ldots,
i_{b+1}$} \\
\end{cases}$$

 There are two remaining
cases: $k$ appears in $\{i_2,\cdots , i_{b+1}\}$ but $k+1$ not, or
$k+1$ appears in $\{i_2,\cdots , i_{b+1}\}$ but $k$ not.  Note for
both cases, the hook distance between $k$ and $k+1$ in the two
hooks $h$ and $\beta$ is the same $\mp k$.  Therefore, the action
of $\sigma_k$ on the respective 2-dimensional subspace is the
same. Since there are $(b-1)$ basis elements $v_i, i\neq k$ in
$\{i_2,\cdots , i_{b+1}\}$, we have a factor of $q^{b-1}$
when comparing to the action of
 $\sigma_k$ on $v_{i_2}\wedge \cdots \wedge v_{i_{b+1}}$.

The second part is proved similarly.  The admissibility condition for
standard Young tableaux reduces the rank by 1. $\hs$

\medskip
In general, Jones-Wenzl sectors $\rho_\lambda^{(k,r)}$ have the
following
properties:
\begin{thm} Let $\lambda$ be an admissible Young diagram which is not of
trivial type.

(1) For each $i$, the image $\rho_{\lambda}^{(k,r)}(\sigma_i)$ has
exactly
two distinct eigenvalues, $-1$ and $q$.

(2) ({\bf Bratteli} diagram) Given a $(k,r)$-admissible Young
diagram $\lambda$ with $n$ boxes, then the restriction of
$\rho_{\lambda}^{(k,r)}$ from $B_n$ to $B_{n-1}$ is the direct sum
of the irreducible representations associated to all
$(k,r)$-admissible Young diagrams $\lambda'$ of size $n-1$
 obtained from $\lambda$ by removing a single corner box.

(3) If $r\ge 5$ and $r\notin \{6,10\}, n\ge 3,$ or $r=10, n\ge 5$,
then the image group of $\rho_{\lambda}^{(k,r)}(B_n)$ is infinite
modulo its center.


\end{thm}

All three statements are in [J2].  The first is obvious from the
construction given above.
One can easily deduce (3) from (1) and (2) given Theorem~\ref{finite}.

\section{Duality of Jones-Wenzl representations}

The Hecke algebra $H_n(q)$ has an automorphism which intertwines
the Jones-Wenzl representations of $H_n(q)$ associated to a pair
of Young diagrams.  This duality was first discovered by F.~Goodman
and H.~Wenzl [GW] and by A.~Kuniba and T.~Nakanishi [KN].
It is called {\it rank-level} duality in conformal
field theory. This duality accounts for the appearance of the
symplectic and orthogonal groups as closed images of certain
Jones-Wenzl representations.

Let $\mathbb N$ denote the set of natural numbers (including $0$).

\begin{definition}
Fix an integer $r>0$.  An \emph{$r$-tile} is a $k\times (r-k)$ matrix
 $T=(t_{ij})_{k\times r-k}$ satisfying the following conditions:

(1) $t_{ij}\in \mathbb N$,

(2) the entries in each row and column are non-increasing,

(3) the difference of any two entries in a single row or column is $\le
1$.
\end{definition}

The relation between $r$-tiles and $(k,r)$-admissible Young
diagrams is given by the following constructions.

\smallskip\noindent
{\bf The $r$-tile $T_{\lambda}$ of a Young diagram $\lambda$:}
Suppose $\lambda=[\lambda_1,\cdots, \lambda_k]$ is a Young diagram
with $k$ rows and $r\geq k+1$.  Let $l=r-k$, and let $T_\lambda$ be the
$k\times l$ matrix with
$$t_{ij}=\left\lfloor\frac{\lambda_i+l-j}{l}\right\rfloor.$$

\smallskip\noindent
{\bf The Young diagram $\lambda_T$ of an $r$-tile $T$:}\ the
 $(k,r)$-admissible Young diagram $\lambda_T$ is a Young
diagram with at most $k$ rows whose $i$th row has
$\sum_{j=1}^{l}t_{ij}$ boxes.

\begin{definition}\par
(1) Given a $(k,r)$-admissible Young diagram $\lambda$, the
\emph{$r$-conjugate of $\lambda$}, denoted $\lambda^{*}_r$, is the
Young diagram associated with the transpose tile of $T_{\lambda}$.

(2) A Young diagram is \emph{$r$-symmetric} if $T_\lambda$ is a
symmetric matrix after discarding all $0$-rows and $0$-columns.

(3) Given a Young tableau $t$ of shape $\lambda$, the
\emph{$r$-conjugate
$t^{*}$} is the tableau of shape $\lambda_r^*$ such that the shape
of $t^{(i)}$ is $r$-conjugate to the shape of ${t^*}^{(i)}$ for all $i$.

\end{definition}

We have the following duality:
\begin{thm}\label{duality}
For any $(k,r)$-admissible Young diagram $\lambda$,
$\rho_{\lambda_{r}^{*}}$ is equivalent to $\chi\otimes
\rho^{*}_{\lambda}$,
where
$\rho^{*}_{\lambda}$ is the contragredient representation of
$\rho_{\lambda}$ and $\chi:B_n
\rightarrow U(1)$ denotes the character with
$\chi(\sigma_i)=-q$.
\end{thm}

{\bf Proof:}  We describe this duality explicitly in terms of
bases. From the definition of the representations $\rho_{\lambda}$
and $\rho_{\lambda_r^{*}}$, the basis elements of the
representation spaces $V_{\lambda}$ and $V_{\lambda_r^{*}}$  are
in 1-1 correspondence by $r$-conjugation of Young tableaux:
$t\leftrightarrow t^{*}$. We define the {\it duality
transformation} $J$ as the linear map $J: V_{\lambda}\rightarrow
V_{\lambda_r^{*}}$ with $J(\vec{v}_{t})=\pm \vec{v}_{t^{*}}$, where
the sign $\pm$ is determined  as follows. Let $t_{0}$ be the
standard vertical tableau of shape $\lambda$.  This is the tableau
in which numbers $1$ through $n$ are filled in one column at a time,
working left to right, and it is not necessarily
admissible.  Each standard tableau $t$ of shape $\lambda$
determines a permutation of $\{1,2,\cdots, n\}$ by comparison to $t_0$.
The sign $\pm$ is the sign of
this permutation.

We show that $\rho_{\lambda_r^{*}}=\chi\otimes \rho^{*}_{\lambda}$
for each braid generator $\sigma_i$.  Given a standard tableau
$t$, there are two cases depending on whether or not $s_i(t)$ is
standard. If
$s_i(t)$ is not standard, then the proof is straightforward.  If
$s_i(t)$ is standard, then
$$\rho_{\lambda}(\sigma_{i})=\left(
\begin{array}{cc}
q-(1+q)\alpha_{t,i} & -\beta_{t,i}\\
-\beta_{t,i}&q-(1+q)(1-\alpha_{t,i})
\end{array} \right).$$
Note that $d_{t^{*},i,i+1}=-d_{t,i,i+1},$ therefore
$\alpha_{t^{*},i}=1-\alpha_{t,i}$.  Since
$\text{det}(\rho_{\lambda}(\sigma_i))=-q$, we have
$$\rho_{\lambda}(\sigma_{i})=(-q)\cdot
\frac{1}{\text{det}(\rho_{\lambda}(\sigma_i))}\cdot
\rho_{\lambda}(\sigma_{i})= \chi\cdot
\rho_{\lambda}^{-1}(\sigma_{i})=\chi \otimes
\rho_{\lambda}^{*}(\sigma_{i}).$$  $\hs$
\begin{corollary}\label{selfdual}\par
(1) If $\lambda$ is $r$-symmetric, then $\dim V_{\lambda}$
is even.

(2) If $\lambda$ is $r$-symmetric, then $\rho_{\lambda}$ is
self-dual up to the character $\chi$. More precisely, suppose
$T=(t_{ij})$ is the $r$-tile of $\lambda$, then if
$\sum_{i>j}t_{ij}$ is odd, $\rho_{\lambda}$ is symplectic up to
$\chi$, and if $\sum_{i>j}t_{ij}$ is even, $\rho_{\lambda}$ is
orthogonal up to $\chi$.

\end{corollary}

{\bf Proof:}  Let us examine more carefully the matrix $J$ representing
the above duality.  First note that $r$-conjugation is an
involution on the basis elements of $V_{\lambda}$ without any
fixed points as long as $\lambda$ has $\geq 2$ boxes.  This
implies (1).  If the sign of $t$ is the same as that of $t^{*}$,
then $J$ is either $\left(
\begin{array}{cc}
0 & 1 \\ 1&0
\end{array} \right)$ or $\left(
\begin{array}{cc}
0 & -1 \\ -1&0
\end{array} \right)$.  Therefore, $J$ defines an orthagonal pairing.
 If the signs of $t$ and $t^{*}$ are different, then $J$
is $\left(
\begin{array}{cc}
0 & 1 \\ -1&0
\end{array} \right)$ or $\left(
\begin{array}{cc}
0 & -1 \\ 1&0
\end{array} \right)$, so $J$ defines a symplectic pairing.
  As $\rho_{\lambda}\cdot
J^{-1}=\chi\otimes {\rho_{\lambda}^{*}}$, up to the character
$\chi$, $\rho$ is either a symplectic or an orthogonal matrix with
respect to either the symplectic form or inner product given by
$J^{-1}$. Checking signs gives (2). $\hs$

The converse of (2) is also true for $r>4$.  This is a slight
refinement of a result of [GW], and we follow the proof given
there.

\begin{thm}
\label{distinct}
Let $r>4$ and $1<k_1,k_2<r-1$.

(1) Let $\lambda_1\in \Lambda_n^{(k_1,r)}$ and $\lambda_2\in
\Lambda_n^{(k_2,r)}$.  If $\lambda_i$ are not of trivial type, then
$\rho_{\lambda_1}$ is equivalent to the tensor
product of $\rho_{\lambda_2}$ with a character of $B_n$ if and only if
$\lambda_1=\lambda_2$.

(2) Let $\lambda_1\in \Lambda_n^{(k_1,r)}$ and $\lambda_2\in
\Lambda_n^{(k_2,r)}$.  If $\lambda_i$ are not of trivial type, then
$\rho_{\lambda_1}$ is equivalent to the tensor
product of $\rho^*_{\lambda_2}$ with a character of $B_n$ if and only if
$\lambda_1=(\lambda_2)_r^*$.
\end{thm}

{\bf Proof.}   For any pair of distinct diagrams $\lambda_1$
and $\lambda_2$, the sets of diagrams of the form $\lambda_1^{(1)}$ and
$\lambda_2^{(1)}$ cannot coincide.  In other words, there exists an
admissible
subdiagram $\mu$ of one of the two, obtained by removing a single box,
which
cannot  be so obtained from the other.  Unless one or both is the Burau
hook
$[n-1,1]$ or its conjugate, $\mu$ is not of trivial type.  If
$\rho_{\lambda_1}$ and $\rho_{\lambda_2}$ are equivalent up to tensoring
by
a character, the same is true of their restrictions to $B_{n-1}$.
We may
therefore proceed by induction, the base case being that in which either
$\lambda_1$ or $\lambda_2$ is $[n-1,1]$ and the other is
$[2,1,\ldots,1]$.
These are not equivalent for $n\ge 4$ by Theorem~\ref{hook}.

Part (2) is an immediate consequence of (1) and
Theorem~\ref{duality}. $\hs$

\section{Closed images of Jones-Wenzl sectors}
In this section, we compute the universal cover $G_1$ of the identity
component $G_0$
of the closure of $\rho_\lambda(B_n)$ for each $\rho_\lambda$ with
infinite
image.  We also give the ambient representation $V$ of $G_0$ (specified
as
a representation of $G_1$.)  Since
$\overline{\rho_\lambda(B_n)}$ is the product of $G_0$ and a group of
scalar
matrices, this is enough information to
determine the actual closure of the image of the sector.

\begin{thm}\label{main}
Fix integers $r$, $n$ such that $r\ge 5$, $r\neq 6$, and $n\ge 3$.  Let
$k$ be
an integer less than $r-1$ and let $\lambda\in\Lambda_n^{(k,r)}$.  We
assume
that $\lambda$ is not of trivial type, and if $r=10$, we assume that
$\lambda$ is neither $[2,1]$ nor $[2,2]$.  Let $G_1$ denote the
universal
cover of the identity component of the closure of $\rho_\lambda(B_n)$
and
$V$, of dimension $N$,
denote the representation space of $\rho_\lambda$ regarded as a
$G_1$-module.  Then

(1) if $\lambda$ is neither $r$-symmetric nor a hook, then $(G_1,V)$
is equivalent to $(SU(N),V_{\varpi_1})$.

(2) if $\lambda$ is a hook with $a+1$ columns and $b+1$ rows, then
$(G_1,V)$ is
equivalent to $(SU(a+b),V_{\varpi_b})$.

(3) if $\lambda$ is not a hook but is $r$-symmetric,
$T_\lambda=(t_{ij})$ is the
$r$-tile of $\lambda$, and $\Sigma=\sum_{i>j}t_{ij}$, then

\quad if $\Sigma$ is even, then $(G_1,V)$ is equivalent to
$(Spin(N),V_{\varpi_1})$;

\quad if $\Sigma$ is odd, then $(G_1,V)$ is equivalent to
$(Sp(N),V_{\varpi_1})$

\end{thm}

The rest of the section is devoted to the proof of this theorem.
We remark that the excluded cases, $r\in\{3,4,6\}$, $r=10$ and
$\lambda\in\{[2,1],[2,2]\}$, or $\lambda$ of trivial type, are
precisely the cases in which the image was already known to be
finite [J2][BW][GJ].

We have already seen that $\rho_\lambda(\sigma_i)$ has two distinct
eigenvalues
whose ratio $-q$ is not $-1$.  Since the braid generators are all
conjugate to
one another, the conjugacy class of $\rho_\lambda(\sigma_i)$
topologically
generates the closure of $\rho_\lambda(B_n)$.  Thus, $G_1$ is simple,
$V$ is
irreducible with highest weight $\varpi$, and $(G_1,\varpi)$ appears on
the list
given in Theorem~\ref{twoeigenvalue}.

\begin{definition}
A pair $(G_1,V)$ consisting of a simply connected simple Lie group and
an
irreducible representation is \emph{standard} if $G_1$ is isomorphic to
$SU(N)$, $Sp(N)$, or $Spin(N)$, and $\dim  V=N$.
\end{definition}

Our main goal is to show that the pairs $(G_1,V)$ arising from diagrams
which are not hooks are standard.  We rule out the other possibilities
offered by Theorem~\ref{twoeigenvalue} by means of two pieces of
information: $\dim  V$, and the closure of $B_{n-1}$ in $G_0$, as
computed by means of the Bratteli diagram.  In order to start the
induction argument, we need to compile results in a number of special
cases.
We begin with hooks.

\begin{proposition}\label{hookcase}
Theorem~\ref{main} holds for all hooks $\lambda$.
\end{proposition}

{\bf Proof:} By Theorem~\ref{hook}, it suffices to consider the
case of Burau hooks $\lambda=[m,1]$.  We use induction on $m$. For
$m=2$ (resp. $m=3$), we can appeal to Theorem~\ref{finite} or to
classical results characterizing all finite subgroups of $GL(2)$
(resp. $GL(3)$) [Ft] to show that $G_0=G_1=SU(2)$ (resp. SU(3))
except when $m=2$ and $r=10$.  For general $m<r$, $\dim
\rho_{[m,1]}=m$, and by the induction hypothesis, $G_0\supset
SU(m-1)$, so $G_0=G_1=SU(m)$. $\hs$

We now consider diagrams $\lambda$ with $\le 7$ boxes which are neither
hooks nor of trivial type.  For $n=4$, $\lambda=[2,2]$, and
$\dim \rho_\lambda=2$, so $G_1=SU(2)$, except when $r=10$, in which
case
$G_1$ is trivial.  For $n=5$, there are two possible diagrams, and
$$\dim \rho_{[3,2]}=\dim \rho_{[2,2,1]}=5,$$
and by Theorem~\ref{twoeigenvalue}, $G_1=SU(5)$ in each case. This
is enough information for the induction argument when $r=5$, so we
now restrict attention to $r\ge7$.  For $n=6$, the diagrams
$[4,2]$, $[3,3]$, $[3,2,1]$, $[2,2,2]$, and $[2,2,1,1]$ give
sectors of dimensions $9$, $5$, $16$, $5$, and $9$ respectively.
Thus, $(G_1,V)$ is obviously standard for each case except the
symmetric diagram $[3,2,1]$, which contains the admissible
subdiagram $[2,2,1]$.  In this case, therefore, $G_0$ contains
$SU(5)$.  It follows that here again, the pair is standard. For
$n=7$, we have $[5,2]$, $[4,3]$, $[4,2,1]$, and $[3,2,2]$ together
with their conjugates; the dimensions are $14$, $14$, $35$, and
$21$ respectively, so Theorem~\ref{twoeigenvalue} implies all are
standard. For $n\ge8$, $\lambda\in\{[4,4],[2,2,2,2]\}$ gives $\dim
\rho_\lambda=14$ and $(G_1,V)$ standard, and otherwise, $\dim
\rho_\lambda>15$.

We can already prove the main theorem in the case that $r=5$.  Indeed,
every $\lambda$ with three rows is $5$-conjugate to one with two, so we
consider only diagrams of the form $[l,m]$, $0\le l-m\le 3$.  By a
Bratteli
diagram computation,
$$\dim \rho_{[l,m]}=\begin{cases} F_{2m-1} &\text{if $l=m$}, \\
F_{2m+1} &\text{if $l=m+1$}, \\
F_{2m+2} &\text{if $m+2\le l\le m+3$,}
\end{cases}$$
where $F_k$ denotes the $k$th Fibonacci number.  If $\dim V=F_{k+1}$
and
$G_0\supset SU(F_k)$, then $G_0=G_1=SU(F_{k+1})$, so the theorem follows
by
induction on $k$.

The general proof of the theorem follows this strategy but is
technically more difficult.  We assume henceforth that $r\ge 7$.

\begin{lemma}\label{oddballs}
The pair $(Spin(8),8)$ never appears among pairs $(G_1,\dim V)$. The
pairs
$(SU(5),10)$, and $(SU(6),15)$ occur only when
$\lambda$ is a hook.
\end{lemma}
{\bf Proof:}
We know already that as $\lambda$ ranges over diagrams which are not
hooks,
$\dim\rho_\lambda$ is never $8$, $10$, or $15$.  When $\lambda$ is a
hook,
$G_1$ is always a special unitary group. $\hs$

\begin{lemma}\label{goursat}
Let $\Lambda\subset\bigcup_k \Lambda_n^{(k,r)}$ denote a set of
diagrams.  Suppose that for
each $\lambda\in\Lambda$, the corresponding pair $(G_1,V)$ is standard.
Let
$\rho_\Lambda$ denote the direct sum of the representations
$\rho_\lambda$,
$\lambda\in\Lambda$.  Then
\begin{equation}\label{onethird}
\mathrm{rank}(\overline{\rho_\Lambda(B_n)}^0)\ge
\frac{\dim\rho_\Lambda}{3}.
\end{equation}
\end{lemma}

{\bf Proof:}
Let $\Lambda'$ denote a maximal subset of $\Lambda$ containing no two
$r$-conjugate
diagrams.
Let $H_\lambda$ denote the quotient of $\overline{\rho_\lambda(B_n)}$ by
its center.  This
is always a simple group, either $PSU(N)$, $PSO(N)$, or $PSp(N)$.  The
closure
of the direct sum $\rho_\lambda\oplus\rho_\mu$ maps to $H_\lambda\times
H_\mu$,
and its
image maps onto each factor.  By Goursat's Lemma, either the image is
the
graph of an isomorphism between $H_\mu$ and $H_\lambda$, or it is the
whole product.
Up to isomorphism, $PSU(N)$ has exactly two non-trivial $N$-dimensional
projective representations, and
they are dual to one another.   By Theorem~\ref{distinct}, if
$\lambda,\mu\in\Lambda'$,
there cannot be an isomorphism $H_\lambda\to H_\mu$ commuting with the
maps from $B_n$,
in the $PSU(N)$ case.  There is only one isomorphism class of
non-trivial projective
$N$-dimensional representations of $PSp(N)$, and the same is true for
$PSO(N)$ when
$N\ge 6$ and
$N\neq 8$.  Thus, again there cannot be an isomorphism $H_\lambda\to
H_\mu$ commuting
with the maps from $B_n$.  By Goursat's lemma, we conclude that the
closure of
$\rho_{\Lambda'}(B_n)$ maps onto $\prod_{\lambda\in\Lambda'}H_\lambda$.
The same is
true a fortiori of the closure of $\rho_\Lambda(B_n)$.  If $\lambda$ is
not
$r$-symmetric, then $H_\lambda$ has rank $N-1\ge 2$, and the sum of the
dimensions of $\rho_\lambda$ and $\rho_{\lambda_r^*}$ is $2N\le 3(N-1)$.
Otherwise the rank of $\rho_\lambda$ is $N/2$ and the contribution of
$\lambda$
to $\dim\rho_\Lambda$ is $N$.  Thus, $\dim\rho_\Lambda$ is at most 3
times
the rank of $\rho_\Lambda(B_n)$. $\hs$

We note that among pairs $(G,V)$ satisfying Theorem~\ref{twoeigenvalue},
the only non-standard ones satisfying
$$\label{ineq}
\mbox{rank}\,G\le \frac{\dim V}{3}$$
are $Spin(7)$ with its spin representation
and $SU(4)$ and $SU(5)$ with their fundamental
representations of dimensions 6 and 10 respectively.  By
Lemma~\ref{oddballs}, these
cases are ruled out for pairs arising from $\rho_\lambda(B_n)$.  We
cannot
proceed immediately by induction, however, since the base cases, which
are
the hooks, do not in general satisfy the inequality (\ref{ineq}).  To
remedy
this, we need to analyze partitions $\lambda$ from which hooks can be
obtained by removing a single box.  We therefore define
$$h_{a,b}=[a+1,\underbrace{1,1,\ldots,1}_b],\ \lambda_{a,b}=
[a+1,2,\underbrace{1,\ldots,1}_{b-1}].$$
Note that the admissibility
of $\lambda_{a,b}$ implies the admissibility of $h_{a,b}$ except in the
case
$a=r-2,b=1$.

\begin{proposition}
If $a+b\ge 5$ and $h_{a,b}$ is admissible, then
$$\dim\rho_{\lambda_{a,b}} \ge (14/5)\dim\rho_{h_{a,b}}.$$
\end{proposition}

{\bf Proof:}
Either $a+b<r-1$ and $h_{a,b}$ has two admissible subdiagrams
with $a+b$ boxes, $h_{a-1,b}$ and $h_{a,b-1}$, or $a+b=r-1$ and there
is only one: $h_{a-1,b}$.
In the first case and if $b>1$,
$\lambda_{a,b}$ has three admissible subdiagrams with
$a+b+1$ boxes, $h_{a,b}$, $\lambda_{a-1,b}$, and $\lambda_{a,b-1}$; in
the
second or if $b=1$, only the first two are admissible.
We proceed by induction, the
proposition being true in the case $a+b=5$ and sharp when $(a,b)=(4,1)$.
Suppose that $n$ is given and the proposition is true when $a+b=n-1$.
Now take
$a+b=n$.  In the first case, if $b>1$,
\begin{eqnarray*}
\dim\rho_{\lambda_{a,b}}&=&\dim\rho_{\lambda_{a-1,b}}
    +\dim\rho_{\lambda_{a,b-1}}+\dim\rho_{h_{a,b}}  \\
&\ge& (14/5)(\dim\rho_{h_{a-1,b}}+
    \dim\rho_{h_{a,b-1}})+\dim\rho_{h_{a,b}}=(19/5)\dim\rho_{h_{a,b}},\\
\end{eqnarray*}
while if $b=1$, than $a\ge4$, so
$$
\dim\rho_{\lambda_{a,1}}=\frac{a^2+3a}{2}\ge\frac{14(a+1)}{5}=\frac{14}{
5}
\dim\rho_{h_{a,1}}.
$$
In the second case,
\begin{eqnarray*}
\dim\rho_{\lambda_{a,b}}&=&\dim\rho_{\lambda_{a-1,b}}
+\dim\rho_{h_{a,b}}\ge (14/5)\dim\rho_{h_{a-1,b}}
+\dim\rho_{h_{a,b}}\\
&=&(19/5)\dim\rho_{h_{a,b}}.\\
\end{eqnarray*}  $\hs$

\begin{proposition}
For any $a,b\ge 1$, $\lambda_{a,b}$ satisfies Theorem~\ref{main}.
\end{proposition}

{\bf Proof:}
By the case analysis following Proposition~\ref{hookcase},
we may take $a+b=n\ge 6$, and we may assume the
proposition is true when $a+b<n$.     The induction hypothesis gives
$\mbox{rank}\,G_1\ge 13$.  Applying Lemma~\ref{goursat} to
$\lambda_{a-1,b}$ and (assuming $b>1$ and $a+b<r-1$) $\lambda_{a,b-1}$,
the induction hypothesis together with Lemma~\ref{onethird} implies that
the
rank of $G_1$ is at least $3/14$ times
the dimension of the representation.  Among the possible pairs $(G_1,V)$
in Theorem~\ref{twoeigenvalue}, only the standard ones satisfy both
conditions. By Lemma~\ref{selfdual}, $G_1$ is
unitary, spin, or orthogonal, depending on which of the conditions in
Theorem~\ref{main} $\lambda_{a,b}$ satisfies.  The proposition follows
by
induction on $n$.
$\hs$

We can now prove Theorem~\ref{main}.

{\bf Proof:} We use induction on $n$.  We may assume that
$\lambda$ is not a hook and that for every admissible tableau with
shape $\lambda$, neither is $\lambda^{(1)}$. Let $\Lambda$ denote
the set of admissible diagrams of the form $\lambda^{(1)}$ for
some admissible tableau.  By inequality (\ref{onethird}),
$$\mbox{rank}\,\overline{\rho_\lambda(B_n)}\ge\mbox{rank}\,
\overline{\rho_\Lambda(B_{n-1})}\ge\frac{\dim\rho_\lambda}{3}.$$
By Lemma~\ref{oddballs}, this inequality together with the fact
that $\lambda$ is not a hook implies that the pair $(G_1,V)$
arising from $\rho_\lambda$ is standard.  The theorem follows by
induction. $\hs$

For completeness,  we point out the closed images of the remaining
cases using Theorem 1.6.  They have all been identified earlier in
[J2][BW][GJ].  As we mentioned earlier, they are all finite
groups. The images for $SU(2), r=4$ are given by Theorem 1.6, (b)
[J2]; $SU(2), r=6$ by Theorem 1.6, (c) cases $(1), (2), (6)$ [BW];
$SU(2), r=10$ and $n=3,4$ by Theorem 1.6, (a) [J2];  The images
for $SU(3), r=6$ are identified first by D.~Goldschmidt and
V.~Jones (see [GJ]), the images are given by Theorem 1.6, (c)
cases $(3), (5)$. The images for $SU(4), r=6$ are the same as
those for $SU(2), r=6$ by rank-level duality.

\section{Distribution of evaluations of Jones polynomials}

In this section, we fix an integer $r\geq 3, r\neq 3,4,6$, and
$q=e^{\pm \frac{2\pi}{r}}$. Given a braid $\sigma \in B_{n}$, let
$\hat{\sigma}$ be the usual closure of $\sigma$.  Then the Jones
polynomial of the link $\hat{\sigma}$ at $q$ is:
$$J(\hat{\sigma},q)=(-1)^{n-1+e(\sigma)}\cdot
q^{-\frac{3e(\sigma)}{2}}\cdot \sum_{\lambda=[\lambda_1,
\lambda_2]\in
\Lambda_{n}^{(2,r)}}\frac{[\lambda_{1}-\lambda_{2}+1]}{[2]}\cdot
Tr(\rho_{\lambda}^{(2,r)}(\sigma)),$$  where $e(\sigma)$ is the
sum of all exponents of standard braid generators appearing in
$\sigma$. In the following,
we denote
$\frac{[\lambda_{1}-\lambda_{2}+1]}{[2]}$
by $w_{\lambda}$.

The sum of exponents $e(\sigma)$ defines a homomorphism from $B_n$
to $\mathbb Z$.  Let $\rho$ denote the direct sum of the
representations $\rho_\lambda$ as $\lambda$ ranges over
$\Lambda_n^{(2,r)}$. Let $G=\overline{\rho(B_n)}\times{\mathbb
Z}_{2r}$. There is a natural map $\rho':\,B_n\to G$ defined by
$\rho'(\sigma)=(\rho(\sigma),r(n-1+e(\sigma))-3e(\sigma)\;\; (mod
\;\; 2r))$. Let
$$T_n:\,\big(\prod_{\lambda\in\Lambda_n^{(2,r)}}U(\lambda)\big)\times
{\mathbb Z}_{2r}\to{\mathbb C}$$ be defined by
$$T_n((u_\lambda),m)=q^{\frac{m}{2}}\sum_{\lambda\in\Lambda_n^{(2,r)}}
w_\lambda Tr(u_\lambda).$$ The definitions are designed so that
$$J(\hat{\sigma},q)=T_n(\rho'(\sigma)).$$ Let $G'\subset G$ denote
the closure of $\rho'(B_n)$.

\begin{lemma}\label{gprime}
If $n\ge 5$, then
$$(G')_0=\prod_{\lambda\in\Lambda_n^{(2,r)}}SU(\lambda),$$ and
$G'=(G')_0 Z(G')$.
\end{lemma}

{\bf Proof:} As $n>4$, a diagram with two rows cannot be
symmetric, nor can two distinct diagrams with two rows be
conjugate to one another.  The computation of $(G')_0$ now follows
immediately from the proof of Lemma~\ref{goursat}.  As $G'$ is a
subgroup of
$$\big(\prod_{\lambda\in\Lambda_n^{(2,r)}}\overline{\rho_\lambda(B_n)}\big)\times
{\mathbb Z}_{2r}$$ and has the same identity component, it
suffices to prove that the latter group is the product of its
identity component and its center.  This is immediate from
Theorem~\ref{simplicity}. $\hs$

\begin{lemma}Let $\mu_{n,k}$ denote probability measure on $\mathbb C$
given
by values of $J(\hat\sigma,q)$, if $\sigma$ is chosen randomly and
uniformly
from (non-reduced) words of length $k$ in the braid generators
$\sigma^{\pm},
\ldots,\sigma_{n-1}^{\pm}\in B_n$.  The weak-* limit of $\mu_{n,k}$ as
$k\to\infty$ is the push-forward of Haar measure on $G'$, ${T_n}_* dg'$.
\end{lemma}

{\bf Proof:}
Let $\nu$ denote the probability measure on $G'$ given by the average of
$\delta$-functions centered at
$\rho'(\sigma_1)^\pm,\ldots,\rho'(\sigma_{n-1})^\pm$.  By [Bh], since
$\rho'(B_n)$ is dense in $G'$, the weak-* limit of the $k$-fold
convolution
$\nu^{*k}$ is Haar measure $dg'$.  Thus the weak-* limit of
${T_n}_*(\nu^{*k})$ is ${T_n}_* dg'$.  $\hs$

The only significance of the choice of the set $\{\sigma_i^\pm\}$ is
that it
generates $B_n$; any other semigroup generators would do as well.
Much more
sophisticated results in ergodic theory can be applied to prove
convergence of measure on more refined ensembles of braids. For
example, the Stein-Nevo theorem [SN] allows the study of reduced
words in the free group.  If $\mu_r$ and $\mu_{r+1}$ are measures
uniformly supported on reduced words in $\gamma_1,\cdots \gamma_m$
and their inverses, then $\frac{1}{2}(\mu_r+\mu_{r+1})$ will also
converge weakly to Haar($G'$).  One may also ask about using the
braid group---not the free group---to count braids and whether a
similar uniformity is obtained.  We do not know at present.

\begin{lemma} If $n\ge r-2$, then
$$\sum_{\lambda\in\Lambda_n^{(2,r)}}w_\lambda^2=
\frac{r}{\sin^2 \frac{2\pi}r}.
$$
\end{lemma}

{\bf Proof:}
There are four cases, depending on the parity of $n$ and $r$.  If both
are even,
the sum in question is
$$[2]^{-2}\sum_{k=0}^{r/2-1}[2k+1]^2=(q-q^{-1})^{-2}\sum_{k=0}^{r/2-1}
\bigl(q^{2k+1}+q^{-1-2k}-2\bigr)
=\frac{r}{\sin^2\frac{2\pi}r}.$$
If $r$ is even and $n$ is odd, the sum is
$$[2]^{-2}\sum_{k=0}^{r/2-2}[2k+2]^2=(q-q^{-1})^{-2}\sum_{k=0}^{r/2-2}
\bigl(q^{2k+2}+q^{-2-2k}-2\bigr)
=\frac{r}{\sin^2\frac{2\pi}r}.$$
If $r$ is odd and $n$ is even, the sum is
$$[2]^{-2}\sum_{k=0}^{r/2-3/2}[2k+2]^2=(q-q^{-1})^{-2}
\sum_{k=0}^{r/2-3/2}
\bigl(q^{2k+2}+q^{-2-2k}-2\bigr)
=\frac{r}{\sin^2\frac{2\pi}r}.$$
Finally, if both are odd,
$$[2]^{-2}\sum_{k=0}^{r/2-3/2}[2k+1]^2=(q-q^{-1})^{-2}
\sum_{k=0}^{r/2-3/2}
\bigl(q^{2k+1}+q^{-1-2k}-2\bigr)
=\frac{r}{\sin^2\frac{2\pi}r}.$$
$\hs$

The fact that $\sum_\lambda w_\lambda^2$ does not depend on the parity
of $n$
has the interesting consequence that the distribution of values of $J$
on
braids of $n$ strands tends to a limit as $n$ goes to $\infty$:

\begin{thm}\label{braid}
The weak-* limit of the sequence of measures ${T_n}_* dg'$ is
the Gaussian distribution
$\frac{1}{2\pi \sigma_{r} }
e^{-\frac{z\bar{z}}{\sigma_{r}}}dzd\bar{z}$, where
$\sigma_r=\frac{r}{\sin^2 2\pi/r}$.
\end{thm}

{\bf Proof:}
By Lemma~\ref{gprime}, we can write $G'=(H\times A)/H\cap A$, where $H$
is
a product of special unitary groups and $A$ is finite and abelian.
Every
representation of $G'$ can be regarded as a representation of $H\times
A$ and
every irreducible representation as an exterior tensor product of an
irreducible
representation of $H$ and an irreducible character of $A$.  In
particular,
the restriction of $T_n$ to $G'$ can be regarded as a function on
$H\times A$:
namely a $w_\lambda$-weighted sum of traces of representations
$\sigma_\lambda\boxtimes\tau_\lambda$, where $\sigma_\lambda$ is the
composition of the standard
representation with the projection onto the factor $SU(\lambda)$ of $H$.

Let $N=\inf_{\lambda\in\Lambda_n^{(2,r)}}\dim\rho_\lambda$.  If
$a_\lambda,
b_\lambda$ are non-negative integers with
$$\sum_{\lambda\in\Lambda_n^{(2,r)}}(a_\lambda+b_\lambda)<N,$$
then
$$\bigotimes_{\lambda\in\Lambda_n^{(2,r)}}
(\sigma_\lambda\boxtimes\tau_\lambda)^{\otimes a_\lambda}\otimes
({(\sigma_\lambda\boxtimes\tau_\lambda)^*}^{\otimes b_\lambda}$$
is isotypic on $Z(H)$ and non-trivial unless $a_\lambda=b_\lambda$ for
all $\lambda$.  In this case, the representation is trivial on $A$, so
the dimension of the space of invariants is
$$\dim\left(
\bigotimes_{\lambda\in\Lambda_n^{(2,r)}}
\sigma_\lambda^{\otimes a_\lambda}\otimes
{\sigma_\lambda^*}^{\otimes a_\lambda}
\right)^H
=\prod_{\lambda\in\Lambda_n^{(2,r)}}
a_\lambda!$$
by the invariant theory of $SU(\lambda)$ [Wl].

Let $\{ X_\lambda \}$ denote a set of independent Gaussian random
variables with distribution $\frac1{2\pi}e^{-z\bar z}dz d\bar z$
indexed by $\lambda\in\Lambda_n^{(2,r)}$.  The expectation is
$$E(X_\lambda^a\bar X_\lambda^b)=
\begin{cases}
a!&\text{if $a=b$}\\
0&\text{otherwise.}\\
\end{cases}
$$ Since all $X_\lambda$, $\lambda\in\Lambda_n^{(2,r)}$, are
independent, if $$X=\sum_{\lambda\in\Lambda_n^{(2,r)}} w_\lambda
X_\lambda,$$ then $$E(X^a\bar X^b)=\int_{G'}
T_n(g')^a\overline{T_n(g')}^b dg' =\int_{\mathbb C} z^a\bar z^b
{T_n}_* dg'$$ whenever $a+b<N$.   As $N$ goes to $\infty$ with
$n$, by [Fe], this implies that each moment of ${T_n}_* dg'$
equals the corresponding moment of the measure $\frac{1}{2\pi
\sigma_{r} } e^{-\frac{z\bar{z}}{\sigma_{r}}}dzd\bar{z}$ of $X$
when $n$ is sufficiently large.  This implies weak convergence by
[Fe] VIII.6 and XV.5. (Actually, the results in [Fe] are stated
only for distributions on $\mathbb R$, but the method works for
${\mathbb R}^n$.)$\hs$

We conclude that if $r$ is a fixed integer $r\ge 5$, $r\neq 6$, then in
the
limit as $n\to\infty$, the distribution of values at $e^{\frac{\pm 2\pi
i}{r}}$
of the Jones polynomial
of a ``random'' link with $n$ strands tends to a fixed Gaussian.  The
variance
of this Gaussian depends on $r$ and grows like $r^3$ as $r\to\infty$.

\begin{thm}
For each $n$ and $k$, let $\mu_{n,k}^{\mathrm{knot}}$ denote the
distribution of values of $J(\hat\sigma,e^{2\pi i/r})$, where
$\sigma$ ranges over those non-reduced words of length $k$ in
$B_n$ for which $\hat\sigma$ is a knot.  If $r=5$ or $r\ge 7$,
then in the weak-* topology,
$$\lim_{n\to\infty}\lim_{k\to\infty}\mu_{n,k}^{\mathrm{knot}}
=\frac1{2\pi\sigma_r}e^{-\frac{z\bar z}{\sigma_r}dz d\bar z}, \
\sigma_r=\frac{r}{\sin^2\frac{2\pi}r}.$$
\end{thm}

{\bf Proof:}
A braid $\sigma$ gives rise to a knot $\hat\sigma$ if and only if the image
of $\sigma$ under the standard quotient map $B_n\to S_n$ is an $n$-cycle.
For each $n\ge 5$ we consider the homomorphism $\phi:\,B_n\to G'\times S_n$
obtained from $\rho'$ and the standard quotient map $B_n\to S_n$.  By Goursat's
lemma, the closure of the image is either all of $G'\times S_n$ or an index-$2$
subgroup.  Applying [Bh] to the topological generators $\phi(\sigma_i^{\pm1})$
of this subgroup, we see that in the large $k$ limit, if we condition on a fixed
element of $S_n$, the resulting distribution on $G'$ approaches one of
three possible limits:
Haar measure $dg'$ on $G'$, twice the restriction of $dg'$ to an index-$2$
subgroup $G'_{\textrm{even}}\subset G'$, or twice the restriction of $dg'$
to the non-trivial coset $G'_{\textrm{odd}}=G'\setminus G'_{\textrm{even}}$.
(Note that the factor of 2 is needed in the last two cases to give a probability
measure.)  The argument of Lemma~\ref{gprime} goes through unchanged when
$G'$ is replaced by $G'_{\textrm{even}}$, so the integral of $z^a\bar z^b$
with respect to ${T_n}_*dg'_{\textrm{even}}$ coincides with the integral
with respect to ${T_n}_*dg'$ when $a+b<N$.  By additivity in measure, the
decomposition
$$dg'=dg'|_{G'_{\textrm{even}}}+dg'|_{G'_{\textrm{odd}}}=
\frac 12dg'_{\textrm{even}}+\frac 12dg'_{\textrm{odd}}$$
gives
$$\int z^a \bar z^b{T_n}_*dg'_{\textrm{odd}}=2\int z^a\bar z^b{T_n}_*dg'
-\int z^a\bar z^b{T_n}_*dg'_{\textrm{even}}=\int z^a\bar z^b{T_n}_*dg'$$
for $a+b<N$.  The theorem now follows from [Fe].
$\hs$

{\bf Remark:}  In [DLL], the evaluations of Jones polynomials at
several roots of unity are plotted for prime knots, or prime
alternating knots up to 13 crossings.  While density still holds
for these cases, we do not know if there exist any limiting
distributions for these ensembles of knots (note that our
filtration in Theorem 5.5 and their filtration for the plotting
are different.)

Another interesting direction is to study subgroups of the braid
groups.  By [Sta], a braid $b$ belonging to $B_{k}(n)$, the $k$-th
stage of the lower central series of the braid group $B_n$,
determines a braid closure $\hat{b}$ whose finite type invariants
vanish through type $k+1$.  Since the groups $SU(m)$ are simple,
if $\rho: B_n\rightarrow SU(m)$ is dense then the restriction
$\rho:B_{k}(n)\rightarrow SU(m)$ is also dense. Thus link
invariants with vanishing invariants of type $\leq k+1$ can
approximate the non-perturbative Jones invariants of an arbitrary
link.  It would be nice to follow this with a uniformity (in
measure) statement,  but this seems to lie outside the scope of
the ergodic theorem we know since in the free group $F_n$,  which
we use to parameterize the braid group, the $k$-th term of the
lower central series $F_k(n)$ is infinitely generated.

Let us now come to the question of the rate of
approximation.  Here to have any kind of general positive answer,
one must restrict to semisimple Lie groups (which fortunately is
where the Jones representations we have studied take their
values). To see this, consider $G=S^1$ and the Liouville number
$\gamma=(\sum_n 10^{-n!})2\pi$, while $\gamma$
generates a dense subgroup and the atomic measure on its partial orbit
converges to the rotationally invariant measure, one must wait an
exceptionally long time for the orbit to come near certain points.
In contrast semisimple groups have a distinctly limited supply of
finite subgroups and nothing similar can occur.  A theorem to
this effect can be found in [Ki][So] and appears in its best form
in [NC].

\begin{thm}
Let $X$ be a set closed under inverse in a compact semisimple Lie
group $G$ (with Killing metrics) such that the group closure
$\langle X \rangle$ is dense in $G$. Let $X_l$ be the words of
length $\leq l$ in $X$, then $X_l$ is an $\epsilon$-net in $G$ for
$l={\mathcal O}(log^2(\frac{1}{\epsilon})),$ i.e.,  for all $g \in
G$, $\text{dist}(g,X_l)<\epsilon$.
\end{thm}

Conjecturally the theorem should still hold for $l={\mathcal
O}(log{(\frac{1}{\epsilon})})$ and there are some number
theoretically special generating sets of $SU(2)$ [GJS] for which
such an estimate for $l$ can in fact be obtained.  Such results
now translate into topological statements:
\begin{corollary}
Given a $\lq\lq$conceivable" value $v$ for the evaluation of Jones
polynomial of $\hat{b}$ at a root of unity, i.e.,  one that lies
in the computed support of the limiting distribution for $b\in
B_n$, the $n$-string braids, to approximate $v$ by $v'$,
$||v-v'||<\epsilon$, it is sufficient to consider braids $b_l'\in
B_n$ of length $l={\mathcal O}(log^2(\frac{1}{\epsilon}))$ with
Jones evaluations $b_l'=v'$, $||v-v'||< \epsilon$.
\end{corollary}

\section{Fibonacci representations}

In this section, we apply the techniques of sections 2 and 4 to
prove a density theorem for a different class of representations.
These arise from Chern-Simons theory for $r=5$ and $G=SO(3)$, what
G.~Kuperberg calls the \emph{Fibonacci TQFT} [KK].

We briefly recall the setup.  The geometric objects we consider
are compact oriented surfaces with boundary, not necessarily
connected, endowed with a parameterization of each boundary
component, i.e., a homeomorphism from $S^1$.  Each boundary
component is labeled with an element of $\{0,2\}$. To each labeled
surface $\Sigma$ there is an associated finite-dimensional Hilbert
space $V_\Sigma$ such that $$V_{\Sigma_1\coprod
\Sigma_2}=V_{\Sigma_1}\otimes V_{\Sigma_2}.$$ If $\Sigma$ is a
labeled surface and $f:\,S^1\to\Sigma$ is a simple closed curve,
we can cut $\Sigma$ along $f(S^1)$.  We call the resulting labeled
surface $\Sigma_{f,a}$ if the two new boundary components are
labeled $a$, and
\begin{equation}\label{sum}
V_\Sigma=V_{\Sigma_{f,0}}\oplus V_{\Sigma_{f,2}}.
\end{equation}
If $Aut(\Sigma)$ denotes the group of isotopy classes of orientation,
label,
and parameterization preserving homeomorphisms $\Sigma\to\Sigma$, there
is
a natural projective unitary action on $V_\Sigma$, provided the Hilbert
space
in question is nonzero.  The restriction of this action to the subgroup
stabilizing the points of $f(S^1)$ decomposes according to equation
(\ref{sum}).  When $\Sigma$ has genus $0$, the projective representation
lifts
canonically to a linear representation.

If $\Sigma$ is a disk with label $a$, then $\dim
V_\Sigma=\delta_{0 a}$. If $\Sigma$ is an annulus with labels $a$
and $b$, then $\dim V_\Sigma=\delta_{a b}$. When $a=b$, it makes
sense to ask for the scalar given by the Dehn twist. If $a=0$, it
is $1$: if $a=2$, it is $\omega=e^{\frac{4\pi i}{5}}$. If $\Sigma$
has genus 0 and 3 boundary components with labels
$a,b,c\in\{0,2\}$, then
\begin{equation}\label{pants}
\dim V_\Sigma=
\begin{cases}
0&\text{if $a+b+c=2$,} \\
1&\text{otherwise.} \\
\end{cases}
\end{equation}

\begin{lemma}
If $\Sigma_{g,m,n}$ has genus $g$ and $m$ (resp. $n$) boundary
components
labelled 0 (resp. 2), then
$$\dim V_{\Sigma_{g,m,n}}=5^{\frac{g-1}{2}}\left\{
\left(\frac{1+\sqrt5}{2}\right)^{g+n-1}+(-1)^{g-1}
\left(\frac{1-\sqrt5}{2}\right)^{g+n-1}\right\}.
$$
\end{lemma}

{\bf Proof:} Immediate by induction. $\hs$

Note that the dimension does not depend on $m$: we can ``cap off''
a boundary component with label 0 by gluing on a disk with label
0.  To simplify bookkeeping, we regard each $V_\Sigma$ as a
projective representation space for $P_{g,m+n}$, the pure mapping
class group for a surface of genus $g$ with $m+n$ boundary
components.  The representation factors through $P_{g,n}$ and is
independent of $m$.  Without abuse of notation, we may therefore
denote it $\rho_{g,n}$.

\begin{thm}\label{fib}
Except when $g+n=1$, $\rho_{g,n}(P_{g,n})$ is dense in $PU(\dim
V_{\Sigma_{g,n}})$.
\end{thm}

The exceptional pairs $(1,0)$ and $(0,1)$ arise in different ways.  In
the first
case, there is a two-dimensional projective representation whose image
is
known to be the icosahedral group; in the second case, there is no
representation
since $V_\Sigma$ is 0-dimensional.  The rest of this section is devoted
to the
proof of the theorem.

\begin{lemma}\label{basecase}
Theorem~\ref{fib} holds for $(g,n)=(0,4)$.
\end{lemma}

{\bf Proof:}  We first compute explicitly the representation of this
case using
[KL].  The representation of a braid generator (in an appropriate basis)
is
 $\left(
\begin{array}{cc}
e^{\frac{4\pi i}{5}} & 0 \\ 0&-e^{\frac{2\pi i}{5}}
\end{array} \right)$.  The fusion matrix is
$\left(
\begin{array}{cc}
\frac{\sqrt{5}-1}{2} & -\sqrt{\frac{\sqrt{5}-1}{2}}
\\ -\sqrt{\frac{\sqrt{5}-1}{2}}&\frac{\sqrt{5}-1}{2}
\end{array} \right)$.
It follows that any finite subgroup of $PU(2)=SO(3)$ can be ruled
out quickly except the icosahedral group.  For this, we compute
the trace of the product of two consecutive braid generators.
This trace cannot arise as the trace of an element of the binary
icosahedral group in the 2-dimensional representation. Therefore,
the image must be dense in $PU(2)$. $\hs$

\begin{proposition}
If $\dim V_{\Sigma_{g,n}}>0$, then $\rho_{g,n}$ is irreducible.
\end{proposition}

{\bf Proof:}
First let $g=0$.  The proposition holds for $n\le 4$.  For $n=5$, we
have a
3-dimensional representation, so it is reducible only if it has an
invariant
line.  Regarding $P_{0,5}$ as a quotient of the braid group $B_5$, we
observe
that $\sigma_1$, $\sigma_2$, and $\sigma_4$ must all fix the line, and
all
three eigenvalues must be the same, either $1$ or $\omega$.  In the
first case,
the line is precisely the subspace of $V_{\Sigma_{0,5}}$ associated to a
loop
with label 0 enclosing the first two boundary components of
$\Sigma_{0,5}$;
it is also the subspace associated to a loop with label 0 enclosing the
last
two boundary components of $\Sigma_{0,5}$.  However, if we cut along
both
loops, we are left with a pair of pants whose labels sum to $2$.  This
is
impossible by (\ref{pants}).  On the other hand, if the eigenvalue is
$\omega$,
the line in question lies in the 2-dimensional space associated to a
loop with
label 2 enclosing the last two boundary components of $\Sigma_{0,5}$,
and this
line is fixed by $\sigma_1$ and $\sigma_2$, contrary to
Lemma~\ref{basecase}.

Now we use induction on $n$.  The dimension of $V_{\Sigma_{0,n}}$ is
$F_{n-1}$, where $F$ denotes the Fibonacci sequence.
We can divide $\Sigma_{0,n}$ by a loop enclosing
the last two boundary components or by a loop enclosing the last three.
In
the first case, we obtain a representation of the loop stabilizer which,
by the induction hypothesis, is a sum of irreducible pieces of
dimensions
$F_{n-2}$ and $F_{n-3}$.  In the second case, we obtain a representation
of
the (different) loop stabilizer which decomposes into irreducible pieces
of
dimension $F_{n-4}$ and $2F_{n-3}$.  As
$$F_{n-4}<F_{n-3}<F_{n-2}<2F_{n-3}, $$
the representation of $P_{0,n}$ is irreducible.

For the higher genus case, we use a similar argument, but in this
case, we choose a non-separating loop and a loop which splits off
a $\Sigma_{1,1}$. In this way, we can write two different
restrictions of $\rho_{g,n}$ as (projectivizations of) a direct
sum of two irreducible representations in two different ways.  The
inequality $$\dim V_{\Sigma_{g-1,n}}<\text{inf} \big(\dim
V_{\Sigma_{g-1,n+1}}, \;\;\; 2\dim V_{\Sigma_{g-1,n}} \big)$$
gives the induction step whenever it holds, which means in every
case except when $g+n \le 3$.  The case $(1,0)$ is well-known.
For $(1,1)$ there is nothing to prove. For $(2,0)$ the
decompositions $5=1+4=2+3$ are different. This leaves the cases
$(1,2)$ and $(3,0)$ which can be handled in the same way as
$(0,5)$ above .$\hs$

We can now prove Theorem~\ref{fib}.  We start with $g=0$ and use
induction.
For $n=5$, Theorem~\ref{twoeigenvalue} implies the desired density.  For
$n\ge6$, $F_{n-2}> \frac{F_{n-1}}{2}$, so any closed subgroup of
$U(F_{n-1})$
acting irreducibly and containing $SU(F_{n-2})$ contains $SU(F_{n-1})$.
Excluding the cases $(1,0)$, $(1,1)$, and $(1,2)$, in each case $g>0$,
$$\dim V_{\Sigma_{g-1,n+2}}>\frac{\dim V_{\Sigma_{g,n}}}{2},$$
so the induction hypothesis together with irreducibility is enough to
give
density.  For $(1,2)$, we use Theorem~\ref{twoeigenvalue}, and there is
nothing
to prove for $(1,0)$ or $(1,1)$.


\begin{thebibliography} {[BK]}

\bibitem[BW]{BW}J. Birman, and B. Wajnryb,
 {\it Markov classes in certain finite quotients of Artin's braid
group,}
  Israel J. Math. 56
(1986), no. 2, 160--178.
\bibitem[Bh]{Bh}R. N. Bhattacharya,
{\it Speed of convergence of the $n$-fold convolution of a probability
measure
on a compact group,}  Z. Wahrscheinlichkeitstheorie und Verw. Gebiete 25
(1972/73), 1--10.
\bibitem[Bl]{Bl} H. Blichfeldt, {\it Finite collineation groups},
Univ. Chicago Press, Chicago, Ill., 1917.
\bibitem[DLL]{DLL}O. Dasbach, T. Le, and X.-S. Lin, {\it Quantum
morphing and the Jones polynomial,} preprint, 2001.
\bibitem[Ft]{Ft} W. Feit, The current situation in the theory of finite
simple groups, {\it Actes du Congr\`es International des
Math\'ematiciens
(Nice, 1970)}, Tome I, 55--93.
\bibitem[Fe]{Fe} W. Feller, An Introduction to Probability Theory and
    its Applications, Volume II, John Wiley \& Sons, New York, 1966.
\bibitem[F]{F}M. H. Freedman, {\it Quantum computation and the
localization of the modular functors}, Foundations of
Computational Mathematics (to appear), quant-ph/0003128.
\bibitem[FKW]{FKW}M. Freedman, A. Kitaev, and Z. Wang, {\it
Simulation of topological field theories by quantum computers},
 Comm. Math. Phys. (to appear), quant-ph/0001071.
\bibitem[FLW]{FLW}M. Freedman, M. Larsen, and Z. Wang, {\it
A modular functor which is universal for quantum computation},
Comm. Math. Phys. (to appear), quant-ph/0001108.
\bibitem[FKLW]{FKLW}M. Freedman, A. Kitaev, M. Larsen, and Z.
Wang, {\it Topological quantum computation}, quant-ph/0101025.
\bibitem[GJ]{GJ}D. Goldschmidt, and V.F.R. Jones, {\it Metaplectic link invariants,}
 Geom. Dedicata 31 (1989), no. 2, 165--191.
\bibitem[GJS]{GJS} A. Gamburd, D. Jakobson, and P. Sarnak,
{\it Spectra of elements in the group ring of ${\rm SU}(2)$.}
 J. Eur. Math. Soc. (JEMS) 1 (1999), no. 1, 51--85.
\bibitem[GW]{GW}F. Goodman, and H. Wenzl, {\it
Littlewood-Richardson coefficients for Hecke algebras at roots of
unity}, Advances in Math., {\bf 82}(1990), 244--265.
\bibitem[Hu]{Hu} J. Humphreys, Introduction to Lie algebras and
    representation theory, Springer-Verlag, New York, 1972.
\bibitem[J1]{J1} V.F.R. Jones, {\it Hecke algebra representations
of braid groups and link polynomial,} Ann. Math., {\bf 126}(1987),
335--388.
\bibitem[J2]{J2}V.F.R. Jones, {\it Braid groups, Hecke algebras and
type $II_1$ factors}, Geometric methods in operator algebras,
Proc. of the US-Japan Seminar, Kyoto, July 1983.
\bibitem[Ki]{Ki}A. Kitaev,  {\it Quantum computations: algorithms
and error correction,} Russian Math. Survey, {\bf 52:61}(1997),
1191--1249.
\bibitem[KL]{KL}L. Kauffman, and S. Lins, {\it Temperley-Lieb Recoupling
theory
and invariants of 3-manifolds}, Ann. Math. Stud., vol. 134.
\bibitem[KK]{KK}A. Kitaev, and G. Kuperberg, work in progress.
\bibitem[KN]{KN}A. Kuniba, and T. Nakanishi, {\it Level-rank duality in
fusion RSOS
models.}
Modern quantum field theory (Bombay, 1990),
344--374, World Sci. Publishing, River Edge, NJ, 1991.
\bibitem[MP]{MP}W. Mckay and J. Patera, Tables of dimensions,
indices, and branching rules for representations of simple Lie
algebras, Lecture notes in pure and applied math., vol 69.
\bibitem[NC]{NC}M. Nielsen and I. Chuang, Quantum Computation and
Quantum
Information, Cambridge Univ. Press, 2000.
\bibitem[Se]{Se}J-P. Serre, {\it Groupes alg\'ebriques associ\'es aux
modules de
Hodge-Tate,}
 Journ\'ees de G\'eom\'etrie Alg\'ebrique de Rennes, Vol.
III, pp. 155--188, Ast\'erisque, 65, Soc. Math. France, Paris, 1979.
\bibitem[So]{So}R. Solvay, private communication.
\bibitem[St]{St}R. Steinberg, Endomorphisms of linear algebraic
groups, Memoir of the AMS, vol. 80.
\bibitem[Sta]{Sta}T. Stanford, {\it Braid commutators and Vassiliev
invariants},
 Pacific J. Math. 174 (1996), no. 1, 269--276
\bibitem[SN]{SN}A. Nevo, and E. Stein, {\it Analogs of Wiener's ergodic
theorems
for semisimple groups. I.}
 Ann. of Math. (2) 145 (1997), no. 3,
565--595.
\bibitem[We]{We}H. Wenzl, {\it Hecke algebras of type $A_n$ and
subfactors}, Invent. Math. {\bf 92}(1988), 349--383.
\bibitem[Wl]{Wl}H. Weyl, The classical groups, Princeton University
Press,
    Princeton, 1939.
\bibitem[Za]{Za}A. Zalesskii, private communication.
\end{thebibliography}
\end{document}